\documentclass[12pt]{article}
\usepackage[english]{babel}
\usepackage{lineno} 
\usepackage{graphicx,multicol}
\usepackage{epic,eepic,epsfig}
\usepackage{caption}

\usepackage{amsmath,amsfonts,amsthm,amssymb}
\usepackage{pstricks,pstricks-add,pst-math,pst-xkey}

\usepackage{algorithm, algorithmic}

\usepackage{graphicx,pstricks,pst-node,amssymb}

\usepackage{tikz}

\pgfdeclarelayer{edgelayer}
\pgfdeclarelayer{nodelayer}
\pgfsetlayers{edgelayer,nodelayer,main}

\tikzstyle{none}=[inner sep=0pt]
\definecolor{hexcolor0xf81e1c}{rgb}{0.973,0.118,0.110}
\definecolor{hexcolor0x3c00ff}{rgb}{0.235,0.000,1.000}

\tikzstyle{whitevertex}=[circle,fill=white,draw=black, scale = 0.5]
\tikzstyle{redvertex}=[circle,fill=hexcolor0xf81e1c,draw=black, scale = 0.5]
\tikzstyle{bluevertex}=[circle,fill=hexcolor0x3c00ff,draw=black, scale = 0.5]
\tikzstyle{greenvertex}=[circle,fill=green,draw=black, scale=0.5]
\tikzstyle{purplevertex}=[circle,fill=magenta,draw=black, scale=0.5]
\tikzstyle{grayvertex}=[circle,fill=white,draw=gray, scale=0.5]
\tikzstyle{blackvertex}=[circle,fill=black,draw=black, scale=0.5]

\tikzstyle{textbox}=[rectangle,fill=none,draw=none]
\tikzstyle{box}=[rectangle,fill=none,draw=black]

\tikzstyle{arc}=[black, ->]
\tikzstyle{grayarc}=[gray, ->]
\tikzstyle{bluearc}=[blue, ->]
\tikzstyle{grayedge}=[draw=gray]
\tikzstyle{blueedge}=[draw=blue]
\tikzstyle{rededge}=[draw=red]
\tikzstyle{edge}=[draw=black]

\tikzstyle{vertex}=[circle, ,fill=white,draw=black, scale=0.5]

\tikzstyle{10circle}=[circle, scale=10.0,draw=black]
\tikzstyle{10oval}=[ellipse, scale=10.0,draw=black]

\usepackage{amsmath,amsfonts,amsthm,amssymb}
\usepackage{pstricks,pstricks-add,pst-math,pst-xkey}

\usepackage{algorithm, algorithmic}

\usepackage{graphicx,pstricks,pst-node,amssymb}

\usepackage{algorithm, algorithmic}

\begin{document}
\date{\today}

\newtheorem{tm}{\hspace{5mm}Theorem}[section]
\newtheorem{prp}[tm]{\hspace{5mm}Proposition}
\newtheorem{dfn}[tm]{\hspace{5mm}Definition}
\newtheorem{lemma}[tm]{\hspace{5mm}Lemma}
\newtheorem{cor}[tm]{\hspace{5mm}Corollary}
\newtheorem{conj}[tm]{\hspace{5mm}Conjecture}
\newtheorem{prob}[tm]{\hspace{5mm}Problem}
\newtheorem{quest}[tm]{\hspace{5mm}Question}
\newtheorem{alg}[tm]{\hspace{5mm}Algorithm}
\newtheorem{sub}[tm]{\hspace{5mm}Algorithm}
\newcommand{\induce}[2]{\mbox{$ #1 \langle #2 \rangle$}}
\newcommand{\2}{\vspace{2mm}}
\newcommand{\dom}{\mbox{$\rightarrow$}}
\newcommand{\ndom}{\mbox{$\not\rightarrow$}}
\newcommand{\compdom}{\mbox{$\Rightarrow$}}
\newcommand{\cdom}{\compdom}
\newcommand{\sdom}{\mbox{$\Rightarrow$}}
\newcommand{\lsd}{locally semicomplete digraph}
\newcommand{\lt}{local tournament}
\newcommand{\la}{\langle}
\newcommand{\ra}{\rangle}
\newcommand{\pf}{{\bf Proof: }}
\newtheorem{claim}{Claim}
\newcommand{\beq}{\begin{equation}}
\newcommand{\eeq}{\end{equation}}
\newcommand{\<}[1]{\left\langle{#1}\right\rangle}

\newcommand{\Z}{\mathbb{$Z$}}
\newcommand{\Q}{\mathbb{$Q$}}
\newcommand{\R}{\mathbb{$R$}}


\title{Comparability and Cocomparability Bigraphs}

\author{Pavol Hell\thanks{School of Computing Science, Simon Fraser University,
Burnaby, B.C., Canada V5A 1S6; pavol@sfu.ca}, 
       Jing Huang\thanks{Department of Mathematics and Statistics, 
      University of Victoria, Victoria, B.C., Canada V8W 2Y2; huangj@uvic.ca}, 
      Jephian C.-H. Lin\thanks{Department of Mathematics and Statistics, 
      University of Victoria, Victoria, B.C., Canada V8W 2Y2; 
      jephianlin@gmail.com},
and Ross M. McConnell\thanks{Computer Science Department, Colorado State 
     University, Fort Collins, CO 80523-1873; rmm@cs.colostate.edu}}
\date{}

\maketitle

\begin{abstract}
We propose bipartite analogues of comparability and cocomparability graphs.
Surprizingly, the two classes coincide. We call these bipartite
graphs cocomparability bigraphs. We characterize cocomparability bigraphs 
in terms of vertex orderings, forbidden substructures, and orientations of their
complements. In particular, we prove that cocomparability bigraphs are precisely
those bipartite graphs that do not have edge-asteroids; this is analogous to Gallai's 
structural characterization of cocomparability graphs by the absence of (vertex-)
asteroids. Our characterizations imply a robust polynomial-time recognition 
algorithm for the class of cocomparability bigraphs. Finally, we also discuss
a natural relation of cocomparability bigraphs to interval containment bigraphs,
resembling a well-known relation of cocomparability graphs to interval graphs.
\end{abstract}


\section{Introduction} \label{intro}

In this paper we propose bipartite analogues of two popular graph classes, namely the comparability and
the cocomparability graphs \cite{golumbic}. Interestingly, the two analogues coincide, and we obtain just
one new class of bigraphs. This class exhibits some features of both comparability and cocomparability
graphs, but the similarities are significantly stronger with the class of cocomparability graphs; therefore
in this paper we call these bigraphs {\em cocomparability bigraphs}. 

Cocomparability graphs are usually defined as the complements of comparability (i.e., transitively orientable) 
graphs, and their name reflects this fact. However, they are a natural and interesting graph class on their own, 
can be defined independently of their complements, and have an elegant forbidden substructure characterization
\cite{golumbic}. Our cocomparability bigraphs bear strong resemblance to these properties.

We remind the reader that many popular graph classes have bipartite analogues. For instance for chordal graphs
there is a well-known bipartite analogue, namely, chordal bipartite graphs, or {\em chordal bigraphs}. They have 
a similar ordering characterization, forbidden substructure characterization, and even a geometric representation 
characterization \cite{gg,huang}.


Interval graph analogues have a more complex history: the bipartite analogues studied first, namely interval bigraphs 
\cite{hkm,muller,sdrw}, do not share many nice properties of interval graphs -- in particular there is no known forbidden 
substructure characterization. A better bipartite analogue of interval graphs turned out to be the interval containment 
bigraphs discussed below. This class has many similar properties and characterizations to the class of interval graphs, 
in particular an ordering characterization, and a forbidden substructure characterization \cite{fhh}.

When considering what constitutes a natural bigraph analogue of a graph class, it turns out best to be guided by the 
ordering characterizations. Especially revealing are the (equivalent) matrix formulation of the ordering characterizations. 
Take the case of chordal graphs. A graph is {\em chordal} if it does not contain an induced cycle of length greater than
three. Chordal graphs are characterized by the existence of a perfect elimination ordering. An ordering $\prec$ of the 
vertices of a graph $G$ is a {\em perfect elimination ordering} if $u \prec v \prec w$ and $uv \in E(G), uw \in E(G)$
implies that $vw \in E(G)$. To consider the matrix formulation, it is most convenient to think of graphs as {\em reflexive}, 
i.e., having a loop at each vertex. (This makes sense, for instance, for the geometric characterization of chordal graphs 
as intersection graphs of subtrees of a tree: since each of the subtrees intersect itself, each vertex has a loop.) The 
matrix condition is expressed in terms of the adjacency matrix of $G$; because of the loops, the matrix has all $1$'s 
on the main diagonal. The {\em $\Gamma$ matrix} is a two-by-two matrix 

\begin{equation}
 \begin{pmatrix} 1 & 1 \\ 1 & 0\end{pmatrix}
\end{equation}

The $\Gamma$ matrix is a {\em principal submatrix} of an adjacency matrix of a reflexive graph if any of the entries
$1$ lies on the main diagonal. Then a perfect elimination ordering of the vertices of $G$ corresponds to a simultaneous
permutation of rows and columns of the adjacency matrix so the resulting form of the matrix has no $\Gamma$ as a
principal submatrix. In other words, {\em a reflexive graph $G$ is chordal if and only if its adjacency matrix can be 
permuted, by simultaneous row and column permutations, to a form that does not have the $\Gamma$ matrix as a 
principal submatrix}. For bipartite graphs, we use the {\em bi-adjacency matrix}, in which rows correspond to vertices 
of one part and columns to vertices of the other part, with an entry $1$ in a row and a column if and only if the two
corresponding vertices are adjacent. Note that this means that re-ordering of the vertices corresponds to independent
permutations of rows and columns. Chordal bigraphs have an ordering characterization \cite{bls} which translates to 
the following matrix formulation. {\em A bipartite graph $G$ is a chordal bigraph if and only if its bi-adjacency matrix 
can be permuted, by independent row and column permutations, to a form that does not have the $\Gamma$ matrix 
as a submatrix.} This indeed yields a class with nice properties analogous to chordal graphs. In particular, a bipartite
graph is a chordal bigraph if and only if it does not contain an induced cycle of length greater than four \cite{gg}.

We next look at the case of interval graphs. A graph is an {\em interval graph} if it is the intersection graph of a family
of intervals in the real line. As for chordality, it is most natural to consider these to be reflexive graphs. Interval graphs
are characterized by the existence of an ordering $\prec$ of $V(G)$ such that if $u \prec v \prec w$ and $uw \in E(G),$
then $vw \in E(G)$. The {\em Slash matrix} is the two-by-two matrix 

\begin{equation}
 \begin{pmatrix} 0 & 1 \\ 1 & 0\end{pmatrix}
\end{equation}

The {\em Slash} matrix is a {\em principal submatrix} of an adjacency matrix of a reflexive graph if either entry $1$ lies 
on the main diagonal. We may now reformulate the ordering characterization as the following matrix characterization.
{\em A reflexive graph $G$ is an interval graph if and only if its adjacency matrix can be permuted, by simultaneous 
row and column permutations, to a form that does not have the $\Gamma$ or the {\em Slash} matrix as a principal 
submatrix}. 

A bipartite graph $G$, with bipartition $(X,Y)$, is an {\em interval containment bigraph} if there is a family of intervals 
$I_v,\ v \in X \cup Y$, such that for $x \in X$ and $y \in Y$, we have $xy \in E(G)$ if and only if $I_x$ contains $I_y$.
The intervals $I_v,\ v \in X \cup Y$, appear to depend on the bipartition $(X,Y)$ in the definition, but it is easy to see
that if $G$ is an interval containment bigraph with respect to one bipartition  $(X,Y)$, it remains so with respect to
any other bipartition \cite{huang}. Simple transformations show \cite{huang} that interval containment bigraphs 
coincide with two other previously investigated classes of bipartite graphs, namely, two-directional orthogonal-ray 
graphs \cite{stu}, and complements of circular arc graphs of clique covering number two \cite{hh}. These graphs
can be characterized by the existence of orderings $\prec_X, \prec_Y$ such that for $u, v \in X$ and $w, z \in Y$, 
if $u \prec_X v, w \prec_Y z$ and $uz, vw \in E(G)$ then $uw \in E(G)$. (Note that $u \prec_X v$ and $w \prec_Y z$ 
mean that the edges $uz, vw$ are ``crossing''.) This implies that {\em a bipartite graph $G$ is an interval containment 
bigraph if and only if its bi-adjacency matrix can be permuted, by independent row and column permutations, 
to a form that does not have the $\Gamma$ or the {\em Slash} matrix as a submatrix.} 

Interval graphs have an elegant structural characterization due to Lekkerkerker and Boland \cite{lb}. An {\em asteroid} 
in a graph is a set of $2k+1$ vertices $v_0, v_1, \dots, v_{2k}$ (with $k \geq 1$) such that for each $i = 0, 1, \dots, 2k$, 
there is a path joining $v_{i+k}$ and $v_{i+k+1}$ whose vertices are not neighbours of $v_i$ (subscript additions are 
modulo $2k+1$). An asteroid with three vertices ($k=1$) is called an {\em asteroidal triple}. The theorem of Lekkerkerker 
and Boland \cite{lb} states  that {\em a graph is an interval graph if and only if it has no induced cycle with more than 
three vertices and no asteroidal triple.} The interval containment bigraphs defined above have an analogous structural 
characterization. As is often the case with bigraph analogues, we must first translate vertex properties into edge properties. 
An {\em edge-asteroid} in a bipartite graph consists of an odd set of edges $e_0, e_1, \dots, e_{2k}$ such that for each 
$i = 0, 1, \dots, 2k$, there is a walk joining $e_{i+k}$ and $e_{i+k+1}$ (including both end vertices of $e_{i+k}$ 
and $e_{i+k+1}$) that contains no vertex adjacent to either end of $e_i$ (subscript additions are modulo $2k+1$).

\begin{figure}[h!]
\begin{center}
\begin{tikzpicture}[>=latex]
                \node [label={}][style=blackvertex] (1) at (0,0) {};
                \node [label={}][style=blackvertex] (2) at (0,-1) {};
                \node [label={}][style=blackvertex] (3) at (-1,.5) {};
                \node [label={}][style=blackvertex] (4) at (-2,.2) {};
                \node [label={}][style=blackvertex] (5) at (1,.5) {};
                \node [label={}][style=blackvertex] (6) at (2,.2) {};
                \node [label={}][style=blackvertex] (7) at (-.7,2) {};
                \node [label={}][style=blackvertex] (8) at (-1.3,2.8) {};
                \node [label={}][style=blackvertex] (9) at (.7,2) {};
                \node [label={}][style=blackvertex] (10) at (1.3,2.8) {};
                \draw [- ] (1) to (2);
                \draw [- ] (4) to (3);
                \draw [- ] (5) to (6);
                \draw [- ] (7) to (8);
                \draw [- ] (9) to (10);
  \draw[-] (-.9,1.5) to [out=90,in=180] (-.35,2.2) to [out=0,in=180] (.35,1.8) 
                  to [out=0,in=-135] (.7,2);`
  \draw[-] (-1,.5) to [out=90,in=180] (-.7,.8) to [out=0,in=180] (.5,-.4) 
                  to [out=0,in=-135] (1.2,.8) ;
\node [style=textbox]  at (-.3, -.5) {$e_0$};
\node [style=textbox]  at (-1.4, 0) {$e_1$};
\node [style=textbox]  at (1.4, 0) {$e_{2k}$};
\node [style=textbox]  at (-1.3, 2.3) {$e_k$};
\node [style=textbox]  at (1.5, 2.3) {$e_{k+1}$};
\end{tikzpicture}
\end{center}
\vspace{-2mm}
\caption{Edge-asteroid}
\label{eas}
\end{figure}
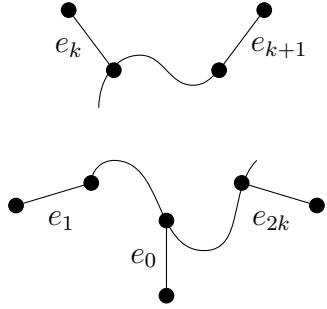


It follows from \cite{fhh} that {\em a bipartite graph $G$ is an interval containment bigraph if and only if it contains no 
induced cycle of length greater than four and no edge-asteroid.}

Armed with these examples we now explore the bipartite analogues of comparability and cocomparability graphs.
It turns out it is most natural to take cocomparability graphs as reflexive, and comparability graphs as irreflexive (i.e., 
without loops). (A hint to the former may be theorem of Gilmore and Hoffman \cite{gh} stating that {\em a graph is an 
interval graph if and only if it is chordal and cocomparability.} Since both interval and chordal graphs are reflexive, it 
makes sense to insist that cocomparability graphs be also reflexive.) The {\em $I_2$ matrix} is the two-by-two identity 
matrix, i.e., 

\begin{equation}
 \begin{pmatrix} 1 & 0 \\ 0 & 1\end{pmatrix}
\end{equation}

Note that the $I_2$ matrix is obtained from the {\em Slash} matrix by the exchange of the two rows, a row permutation. 
The $I_2$ matrix is a {\em principal submatrix} of an adjacency matrix of an irreflexive graph if either of the entries 
$0$ lies on the main diagonal. It is well known that a graph $G$ is a comparability graph if and only if it has an ordering 
$\prec$ of $V(G)$ such that $u \prec v \prec w$ and $uv \in E(G), vw \in E(G)$ implies that $uw \in E(G)$ \cite{skrien}. 
This has a natural matrix formulation as follows. {\em An irreflexive graph is a comparability graph if and only if its adjacency 
matrix can be permuted, by simultaneous row and column permutations, to a form that does not have the $I_2$ matrix as 
a principal submatrix.} Correspondingly, we define a bipartite graph $G$ to be a {\em comparability bigraph} if its bi-adjacency 
matrix can be permuted, by independent row and column permutations, to a form that does not have the $I_2$ matrix as a 
submatrix. 

The situation is similar for cocomparability graphs. Obviously, a graph is a cocomparability graph if and only if it has an 
ordering $\prec$ of $V(G)$ (called a {\em cocomparability ordering}) such that $u \prec v \prec w$ and $uv \not\in E(G),$
$vw \not\in E(G)$ implies that $uw \not\in E(G)$. 
In matrix terms, {\em a reflexive graph is a cocomparability graph if and only if its adjacency matrix can be permuted, by 
simultaneous row and column permutations, to a form that does not have the {\em Slash} matrix as a principal submatrix.} 
Thus we define a bipartite graph $G$ to be a {\em cocomparability bigraph} if its bi-adjacency matrix can be permuted, by 
independent row and column permutations, to a form that does not have the {\em Slash} matrix as a submatrix. By the 
previous observation about the relation of the $I_2$ and the {\em Slash} matrices, we see that by reversing the rows of 
a bi-adjacency matrix without the $I_2$ matrix, we obtain a matrix without the {\em Slash} matrix. Hence we can conclude 
the following fact mentioned earlier.

\begin{prp}\label{first}
A bipartite graph is a comparability bigraph if and only if it is a cocomparability bigraph.
\qed
\end{prp}

Consequently, we shall call our class just by one name. We choose to call these bipartite graphs cocomparability 
bigraphs, as they exhibit more similarities with cocomparability graphs. 

\begin{figure}[h!]
\begin{center}
\begin{tikzpicture}[>=latex]
                \node [label={below:$y'$}][style=blackvertex] (1) at (0,0) {};
                \node [label={above:$x$}][style=blackvertex] (2) at (0,1.5) {};
                \node [label={below:$y$}] [style=blackvertex] (3) at (1.5,0) {};
                \node [label={above:$x'$}][style=blackvertex] (4) at (1.5,1.5) {};
                \draw [- ] (1) to (4);
                \draw [- ] (2) to (3);
                \draw[dotted] (1) -- (2);
                \draw[dotted] (3) -- (4);
\end{tikzpicture}
\end{center}
\vspace{-2mm}
\caption{The forbiddden subgraph $S$ in $G$.}
\text{(The top and bottom vertices are ordered left-to-right according to $\prec_X, \prec_Y$ respectively.)}
\label{slash}
\end{figure}

For future reference, we reformulate the matrix definition in terms of vertex ordering. The {\em Slash} matrix corresponds 
to the pattern $S$ in Figure~\ref{slash}, in the following sense: in the figure we have a bipartition of a bipartite graph
into vertices in the upper row, coming from a set $X$, and those in the lower row, coming from a set $Y$. The set $X$ 
will correspond to the rows of the matrix and the set $Y$ to the columns of the matrix. The independent orderings of the
rows and columns yield two orderings, the ordering $\prec_X$ of the set $X$, and the ordering $\prec_Y$ of the set
$Y$. The depicted pattern $S$ describes precisely the presence of a {\em Slash} submatrix. Therefore, for a bipartite 
graph $G$, with bipartition $(X,Y)$, we say that the pair of orderings $\prec_X$, $\prec_Y$ is {\em S-free}, if for all 
$u, v \in X$ and $w, z \in Y$, if $u \prec_X v, w \prec_Y z$ and $uz, vw \in E(G)$ then at least one of $uw, vz$ is in 
$E(G)$. (We emphasize that being $S$-free is a property of the {\bf pair} $\prec_X$, $\prec_Y$.) Therefore, {\em a 
bipartite graph, with bipartition $(X,Y)$, is a cocomparability bigraph if and only if there exists an $S$-free pair of
orderings $\prec_X, \prec_Y$.} Moreover, it is easy to see that a bipartite graph with a bipartition $(X,Y)$ has an
$S$-free pair of orderings $\prec_X, \prec_Y$ if and only if this is true in any other bipartition.

Suppose $G$ is a bipartite graph with bipartition $(X,Y)$. We define two edges $xy, x'y'$ with $x, x' \in X, y, y' \in Y$ 
to be {\em independent} if they are disjoint and neither $xy'$ nor $x'y$ is an edge of $G$. Given a pair of orderings 
$\prec_X$, $\prec_Y$ of $X, Y$ respectively, we say that two edges $xy, x'y'$ with $x, x' \in X, y, y' \in Y$ {\em cross} 
if $x \prec_X x'$ and $y' \prec_Y y$, or $x' \prec_X x$ and $y \prec_Y y'$. It is clear that the pair $\prec_X$, $\prec_Y$ 
is $S$-free if and only if no two independent edges cross.

We define the {\em independence graph} $I(G)$ of $G$, whose vertices are the edges of $G$, and two edges are 
adjacent just if they are independent. Thus the complement $\overline{I(G)}$ has two edges of $G$ adjacent if and 
only if they share an end or are joined by at least one other edge. Therefore an edge-asteroid in $G$ corresponds 
to a set of vertices and paths in $\overline{I(G)}$, and conversely.

\begin{prp} \label{stars}
Let $G$ be a bipartite graph. Then $G$ contains an edge-asteroid if and only if $\overline{I(G)}$ 
contains an asteroid. \qed
\end{prp}

Let us now point out some similarities of cocomparability bigraphs and cocomparability graphs.

Cocomparability graphs admit an elegant forbidden structure characterization in terms of asteroids. 
{\em A graph is a cocomparability graph if and only if it has no asteroids} \cite{gallai}. Note that it can 
be deduced from this that {\em a graph is an interval graph if and only if it is both a chordal graph 
and a cocomparability graph} \cite{gh}.

Our main theorem in this paper asserts that {\em a bipartite graph is a cocomparability 
bigraph if and only if it does not contain an edge-asteroid.} It implies that {\em a bipartite
graph is an interval containment bigraph if and only if it is both a chordal bigraph and 
a cocomparability bigraph.} These two results strongly resemble the corresponding 
statements for cocomparability graphs and interval graphs discussed above, and make
a good case that these are indeed the right analogues.

This paper is organized as follows. In Section~\ref{forbs} we translate the ordering
properties into properties of orientations of the complements, and we identify two 
forbidden structures for the existence of such orientations, namely, invertible pairs
and edge-asteroids. We show that each of these can be used to characterize bigraphs 
whose complements have suitable orientations. In Section~\ref{characs} we prove 
that if a bigraph contains no invertible pair then it is a cocomparability bigraph. It follows 
that each of the forbidden structures also characterizes cocomparability bigraphs. 
Finally, in Section~\ref{conseq}, we show that cocomparability bigraphs are recognizable 
in polynomial time and point out some consequences of our characterizations of 
cocomparability bigraphs.

\section{Orientations and their obstructions} \label{forbs}

Let $G$ be a bipartite graph, with bipartition $(X,Y)$. Recall that $G$ is a cocomparability 
bigraph if it has an $S$-free pair of orderings $\prec_X$, $\prec_Y$ of $X, Y$ respectively. 
(In other words, orderings in which no two independent edges cross.) The {\em bipartite complement}
$G'$ of $G$ has the same vertices as $G$, and the same bipartition $(X,Y)$, and $xy, x \in X,$ 
$y \in Y$, is an edge of $G'$ if and only if it is not an edge of $G$. Note that it follows from 
Proposition \ref{first} that $G$ is a cocomparability bigraph if and only if $G'$ is one.

We now consider the ordinary complement $\overline{G}$. Note that both $X$ and $Y$ are complete 
subgraphs of $\overline{G}$, and the edges in $\overline{G}$ between $X$ and $Y$ are 
precisely the edges of the bipartite complement $G'$. Denote by ${\cal F}$ the set of pairs
$(a,b)$ where $a \neq b$ and both $a, b$ are in $X$ or both are in $Y$. We will not distinguish
between the pair $(a,b)$ and the corresponding edge $ab$ of $\overline{G}$.
in the two complete subgraphs induced by $X$ and $Y$. We will consider mixed graphs $\vec{G}$,
that are partial orientations of $\overline{G}$ in which the edges of $G'$ remain undirected, 
and each edge of ${\cal F}$ is given an orientation. For simplicity, we will refer to such mixed 
graphs as orientations of $\overline{G}$. Thus from now on the words {\em orientation of}
$\overline{G}$ will mean a mixed graph in which the edges of $\overline{G}$ between $X$ 
and $Y$ are undirected and the edges of $\overline{G}$ inside $X$ and $Y$ are oriented. 
The words {\em partial orientation of} $\overline{G}$ will mean a mixed graph in which all the 
oriented edges are inside ${\cal F}$. (In other words, some edges of ${\cal F}$ may be left
undirected.)

Let $xy, x'y'$, with $x, x' \in X, y, y' \in Y$ be two edges of $\overline{G}$ 
that are independent in the bipartite complement $G'$. (I.e., neither $xy'$ nor $x'y$ are edges of 
$G'$.) We say that the edges $xy, x'y'$ {\em agree} in $\vec{G}$ if both $xx', yy'$ are in $\vec{G}$ 
or both $x'x, y'y$ are in $\vec{G}$. We say that the mixed graph $\vec{G}$ is a {\em $T$-free 
orientation} of $\overline{G}$ if no two independent edges of $G'$ agree in $\vec{G}$. In other 
words, $\vec{G}$ is a $T$-free orientation of $\overline{G}$ if and only if it does not contain as 
induced subgraph the mixed graph $T$ depicted in Figure~\ref{2-trans-forb}.

\vspace{-4mm}
\begin{center}
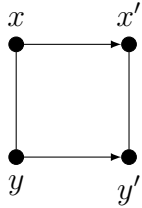
\begin{figure}[htb]
\center
\begin{tikzpicture}[>=latex]
\begin{pgfonlayer}{nodelayer}
                \node [label={above:$x$}][style=blackvertex] (1) at (0,1.5) {};
                \node [label={below:$y$}][style=blackvertex] (2) at (0,0) {};
                \node [label={above:$x'$}][style=blackvertex] (3) at (1.5,1.5) {};
                \node [label={below:$y'$}][style=blackvertex] (4) at (1.5,0) {};
                \draw[ ] (1) to (2);
                \draw[ ] (3) to (4);
                \draw[style=arc] (1) to (3);
                \draw[style=arc] (2) to (4);
\end{pgfonlayer}
\end{tikzpicture}
\caption{\label{2-trans-forb}The forbidden mixed subgraph $T$ in $\vec{G}$}
\end{figure}
\end{center}

We note that $T$ is essentially the complement of $S$.
Suppose $\prec_X$, $\prec_Y$ is an $S$-free pair of orderings of $X, Y$ respectively. We can 
orient the edges of the complement $\overline{G}$ as follows: if $x \prec_X x'$, then the edge 
$xx'$ is oriented from $x$ to $x'$, and similarly for edges $yy'$ in $Y$. It is clear that the resulting 
orientation $\vec{G}$ is $T$-free; moreover, it is also {\em acyclic}, i.e., there is no directed cycle 
in $X$ or in $Y$. Conversely, if we have an acyclic $T$-free orientation $\vec{G}$ of $\overline{G}$, 
any topological ordering (linear extension) of the orientations in $X$ and in $Y$ will give a pair of 
orderings $\prec_X, \prec_Y$ that is $S$-free.

We draw the following conclusion.

\begin{prp} \label{2-tr}
Let $G$ be a bipartite graph with bipartition $(X,Y)$. Then the following statements are equivalent:

\begin{enumerate}
\item
$G$ is a cocomparability bigraph 
\item
$\overline{G}$ has an acyclic $T$-free orientation. 
\end{enumerate}
\end{prp}

We will first study when $\overline{G}$ has any $T$-free orientation. There are two natural obstructions
for this to happen.

We say that two walks in $G$ that both begin in $X$ or both begin in $Y$ are {\em congruent} 
if they have the same length, and if for each $i$ their $i$-th edges are independent. A walk in $G$
is an $(a,b)$-{\em walk} it it starts in $a$ and ends in $b$. A pair of vertices $u, v$ in $G$ is called 
an {\em invertible pair} if there exist congruent walks $W, W'$ where $W$ is a $(u,v)$-walk and $W'$ 
is a $(v,u)$-walk.

\begin{prp} \label{invertible_pair}
If a bipartite graph $G$ contains an invertible pair, then $\overline{G}$ does not have a $T$-free 
orientation.
\end{prp}

\pf Suppose, for a contradiction, that $u, v$ is an invertible pair, with walks $W, W'$ in $G$ in
which consecutive pairs of edges are independent, and that $\vec{G}$ is a $T$-free orientation
of $\overline{G}$. Assume without loss of generality that the first vertex $u$ of $W$ is oriented 
towards the first vertex $v$ of $W'$ in $\vec{G}$. Then it can be seen by induction on $i$ that the 
$i$-th vertex of $W$ must be oriented  the $i$-th vertex of $W'$ in $\vec{G}$. Indeed, suppose 
that $u_iu_{i+1}$ and $v_iv_{i+1}$ are independent edges in $G$ and that $u_iv_i$ and 
$v_{i+1}u_{i+1}$ are oriented edges of $\vec{G}$. Since $u_i,v_{i+1}, u_{i+1}v_i$ are 
undirected edges of $\vec{G}$, we have two independent edges of $G'$ that agree in 
$\vec{G}$, i.e., the orientation $\vec{G}$ is not $T$-free, a contradiction. We conclude that
the last vertex $v$ of $W$ must be oriented towards the last vertex $u$ of $W'$, which is the
opposite of what we started with. \qed

\begin{cor} \label{invertible_p}
If a bipartite graph $G$ contains an invertible pair, then $G$ is not a cocomparability bigraph.
\end{cor}

\begin{prp} \label{Easteroid}
If a bipartite graph $G$ contains an edge-asteroid, then it contains an invertible pair.
\end{prp}

\pf Suppose that $x_0y_0, x_1y_1, \dots, x_{2k}y_{2k}$ form an edge-asteroid in $G$ 
where $x_i\in X$ and $y_i\in Y$ for all $i$, together with walks joining $x_{i+k}y_{i+k}$ 
and $x_{i+k+1}y_{i+k+1}$ that contain $x_{i+k}, y_{i+k}, x_{i+k+1}, y_{i+k+1}$ but no 
vertex adjacent to either of $x_i$ and $y_i$. Let the walk between $x_{i+k}y_{i+k}$ 
and $x_{i+k+1}y_{i+k+1}$ be $v_1v_2 \dots v_k$, where $v_1 = x_{i+k}, v_2 = y_{i+k}$, 
$v_{k-1} = y_{i+k+1}$, and $v_k = x_{i+k+1}$. Then consider the walk $u_1u_2 \dots u_k$, 
with $u_j = x_i$ and $u_{j+1} = y_i$ for each odd $j$. Since $u_jv_{j+1}, v_ju_{j+1} \notin E(G)$ 
for each $1 \leq j \leq k-1$, the two walks $u_1u_2 \dots u_k$ and $v_1v_2 \dots v_k$ are 
congruent. That is, for each $0 \leq i \leq 2k$, there exist congruent walks $W_{i+k}$ and
$W'_i$, where the former walk goes from $x_{i+k}$ to $x_{i+k+1})$, and the latter walk
goes from $x_i$ to $x_i$. By concatenating the walks $W_0', W_0, W_1', W_1, \dots, W_{k-1}, W_k'$, 
we obtain a walk from $(x_0$ to $x_k$, and by concatenating the walks 
$W_k, W'_{k+1}, W_{k+1}, \dots, W'_{2k}, W_{2k}$ we obtain a walk from $x_k$ to $x_0$;
these two walks are congruent, and thus $x_0, x_k$ is an invertible pair in $G$.
\qed

We observe for future reference that the invertible pairs constructed above remain invertible pairs
even if the edges $x_0y_0, x_1y_1, \dots, x_{2k}y_{2k}$ in the edge asteroid are not required to
be distinct, as long as there are walks joining $x_{i+k}y_{i+k}$ and $x_{i+k+1}y_{i+k+1}$ that 
contain $x_{i+k}, y_{i+k}, x_{i+k+1}, y_{i+k+1}$ but no vertex adjacent to either of $x_i$ and $y_i$.
We call such a set of edges a {\em weak edge-asteroid}.

\begin{prp} \label{almostgood}
If $I(G)$ is a comparability graph, then $\overline{G}$ has a $T$-free orientation.
\end{prp}

\pf Let $(X,Y)$ be a bipartition of $G$, and let $\prec$ be a transitive orientation of $I(G)$. 
We orient $\overline{G}$ as follows. Suppose $xy, x'y'$ are two independent edges of $G$;
thus they are adjacent in $I(G)$. Suppose $xy \prec x'y'$ in the transitive orientation of $I(G)$.
Then we put $xx' \in \vec{G}$ and $y'y \in \vec{G}$. Note that this will not create a copy of $T$ 
on $x, x', y, y'$ because in $\vec{G}$ we have the directed edges $xx', yy'$ and the undirected 
edges $xy', x'y$. Any remaining undirected pairs $xx', yy'$ may be oriented arbitrarily. It remains to
show that this is an orientation, i.e., that no edge of $\overline{G}$ inside $X$ or inside $Y$ has
been oriented in both directions. Without loss of generality suppose that this happened for an
edge $xx'$ inside $X$; it was oriented from $x$ to $x'$ because $xy \prec x'y'$, and oriented
from $x'$ to $x$ because $x'z' \prec xz$, in $I(G)$. Note that all of $xy', xz', x'z, x'y$ are non-edges
of $G$, since $xy, x'y'$ and $xz, x'z'$ are independent pairs of edges. This implies that $I(G)$ also
contains edges between $xy$ and $x'z'$ and between $xz$ and $x'y'$. However, $I(G)$ does not 
contain edges between $xy$ and $xz$, and between $x'y'$ and $x'z'$, as those pairs intersect and 
hence are not independent. This contradicts the transitivity of $\prec$. If the edge between $xy$ and 
$x'z'$ has $xy \prec x'z'$ then by transitivity we would have $xy \prec x'z' \prec xz$, contradicting
the fact that there is no edge between $xy$ and $xz$. If it has $x'z' \prec xy$, then $x'z' \prec xy \prec x'z'$,
also a contradiction.
\qed.

Combining Propositions \ref{invertible_pair}, \ref{Easteroid} and \ref{almostgood} 
we obtain the following equivalences. These statements verify that (i) implies (ii), (ii) implies
(iii), (iii) implies (iv), and obviously (iv) implies (v). Proposition \ref{stars} shows that (v) and
(vi) are equivalent, and the equivalence of (vi) and (i) follows by Gallai' theorem \cite{gallai}.

\begin{tm} \label{charac}
The following statements are equivalent for a bipartite graph $G$:
\begin{description}
\item{(i)} $I(G)$ is a comparability graph;
\item{(ii)} $\overline{G}$ has a $T$-free orientation;
\item{(iii)} $G$ does not contain an invertible pair;
\item{(iv)} $G$ does not contain a weak edge-asteroid;
\item{(v)} $G$ does not contain an edge-asteroid;
\item{(vi)} $\overline{I(G)}$ is does not contain an asteroid.
\qed
\end{description}
\end{tm}

We note that the statement $2$ of Proposition \ref{stars} clearly implies the statement $(ii)$ of 
Theorem \ref{charac}, thus any the properties in Proposition \ref{2-tr} imply all the properties in 
Theorem \ref{charac}. In the next section we prove the converse, showing that all these properties 
are in fact equivalent.

\section{Acyclic orientations} \label{characs}

We will show in this section that a bipartite graph satisfying any of the equivalent
conditions in Theorem \ref{charac} is a cocomparability bigraph. (The converse 
has been observed at the end of the previous section.) To do this, we will show 
that the absence of invertible pairs in $G$ allows us to construct an acyclic $T$-free 
orientation of $\overline{I(G)}$.

Let $G$ be a bigraph with bipartition $(X,Y)$.

For $(a,b), (f,g) \in \cal F$, we say that $(a,b)$ {\em implies} $(f,g)$,
and write $(a,b) \Gamma (f,g)$, if there exist congruent $(a,f)$- and $(b,g)$-walks. 
It is easy to verify that $\Gamma$ is an equivalence relation on $\cal F$. 
An equivalence class of this relation will be called an {\em implication class}.
It follows from this definition that there is an implication class contains both $(a,b)$ 
and $(b,a)$ if and only if $a, b$ is an invertible pair. Note that $(a,b) \Gamma (f,g)$ 
if and only if $(b,a) \Gamma (g,f)$. 

Two congruent walks \[a_1a_2 \dots a_{k-1}a_k \ \mbox{and}\ b_1b_2 \dots b_{k-1}b_k\]
are called {\em standard} if for each $i=1, 2, \dots, k-2$ we have $a_i=a_{i+2}$ or $b_i=b_{i+2}$.
It is easy to see that if there exist congruent $(a,f)$- and $(b,g)$-walks, then there exist standard
congruent $(a,f)$- and $(b,g)$-walks. Indeed, suppose that for some $i=1, 2, \dots, k-2$, we
have $a_i \neq a_{i+2}$ and $b_i \neq b_{i+2}$. Note that we must have $a_i \neq b_{i+2}$,
since $a_i$ is not adjacent to $b_{i+1}$ in $G$ but $b_{i+2}$ is. Similarly, $b_i \neq a_{i+2}$.
So $a_i, a_{i+2}, b_i, b_{i+2}$ are all distinct. Thus the following two walks are congruent:
\[a_1a_2 \dots a_ia_{i+1}a_ia_{i+1}a_{i+2} \dots a_{k-1}a_k\ \mbox{and}\ 
b_1b_2 \dots b_ib_{i+1}b_{i+2}b_{i+1}b_{i+2} \dots b_{k-1}b_k.\] 
Continuing this way, we make sure the walks are standard.

For the remainder of this section, let $G$ be a bigraph with bipartition $(X,Y)$ and
\[a_1a_2 \dots a_{k-1}a_k\ \mbox{and}\ b_1b_2 \dots b_{k-1}b_k\]
congruent walks, with an odd $k$, that begin in $X$ and end in $Y$.
Let $C_X$ denote the set of all vertices $a_i$ and $b_i$ with odd $k$,
i.e., all vertices of both paths that lie in $X$, and let $C_Y$ be defined
analogously as the set of all vertices of both paths that lie in $Y$.

\begin{prp} \label{prelem}
Suppose $c_1 \in X, c_2 \in Y$ are two vertices satisfying either of the 
following two conditions:
\begin{itemize}
\item $c_1c_2 \in E(G)$, but $c_1$ is not adjacent to any vertex in $C_Y$ 
and $c_2$ is not adjacent to any vertex in $C_X$;
\item $c_1c_2 \notin E(G)$, but $c_1$ is adjacent to every vertex in $C_Y$ 
and $c_2$ is adjacent to every vertex in $C_X$.
\end{itemize}

Then $(a_1,c_1) \Gamma (a_k,c_1)$ and $(c_1,b_1) \Gamma (c_1,b_k)$.
\end{prp}

\pf Suppose the first condition holds. Then 
\[a_1a_2a_3a_4 \dots a_{k-1}a_k\ \mbox{and}\ c_1c_2c_1c_2 \dots c_2c_1\]
are congruent walks, thus $(a_1,c_1) \Gamma (a_k,c_1)$. 
Similarly, 
\[c_1c_2c_1c_2 \dots c_2 c_1\ \mbox{and}\ b_1b_2b_3b_4 \dots b_{k-1}b_k\]
are congruent walks, thus $(c_1,b_1) \Gamma (c_1,b_k)$.

Suppose now that the second condition holds. Then 
\[a_1c_2a_3c_2 \dots c_2a_k\ \mbox{and}\ c_1b_2c_1b_4 \dots b_{k-1}c_1\]
are congruent walks, thus $(a_1,c_1) \Gamma (a_k,c_1)$. 
Similarly, 
\[c_1a_2c_1a_4 \dots a_{k-1}c_1\ \mbox{and}\ b_1c_2b_3c_2 \dots c_2b_k\]
are congruent walks, thus $(c_1,b_1) \Gamma (c_1,b_k)$.
\qed

We call a pair $(u,v) \in \cal F$ {\em relevant} if its implication class contains at least 
two pairs. Note that $(u,v)$ is relevant if and only if $(v,u)$ is relevant.

\begin{figure}[h]
\begin{center}
\begin{tikzpicture}
\node [style=blackvertex] [label={left:$a_1$}] (a1) at (0,2) {};
\node [style=blackvertex] [label={right:$b_1$}] (b1) at (2,2) {};
\node [style=blackvertex] [label={$c_1$}] (c1) at (1,1.5) {};
\node [style=blackvertex] [label={left:$a_2$}] (a2) at (0,0) {};
\node [style=blackvertex] [label={right:$b_2$}] (b2) at (2,0) {};
\node [style=blackvertex] [label={below:$c_2$}] (c2) at (1,-0.5) {};
\draw (a1)--(a2);
\draw (b1)--(b2);
\draw (c1)--(c2);
\draw[dotted] (a1)--(c2)--(b1)--(a2)--(c1)--(b2)--(a1);
\end{tikzpicture}
\hfil
\begin{tikzpicture}
\node [style=blackvertex] [label={left:$a_1$}] (a1) at (0,2) {};
\node [style=blackvertex] [label={right:$b_1$}] (b1) at (2,2) {};
\node [style=blackvertex] [label={$c_1$}] (c1) at (1,1.5) {};
\node [style=blackvertex] [label={left:$a_2$}] (a2) at (0,0) {};
\node [style=blackvertex] [label={right:$b_2$}] (b2) at (2,0) {};
\node [style=blackvertex] [label={below:$c_2$}] (c2) at (1,-0.5) {};
\draw (a1)--(c2)--(b1)--(b2)--(c1)--(a2)--(a1);
\draw[dotted] (a1)--(b2);
\draw[dotted] (b1)--(a2);
\draw[dotted] (c1)--(c2);
\end{tikzpicture}
\end{center}
\vspace{-2mm}
\caption{An illustration of the proof of Proposition~\ref{mainlem}}
\label{figmainlem}
\end{figure}
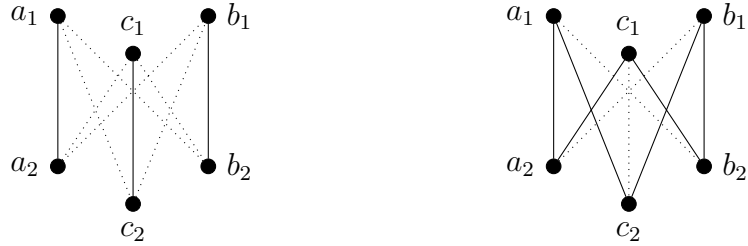

\begin{prp} \label{mainlem}
Suppose in addition that the above walks
\[a_1a_2 \dots a_{k-1}a_k\ \mbox{and}\ b_1b_2 \dots b_{k-1}b_k\]
are standard. 

Suppose $c_1 \in X$ satisfies the following properties:
\begin{itemize}
\item
the pair $(a_1,c_1)$ is relevant;
\item
$(a_1,c_1) \not\Gamma (a_1,b_1)$;
\item
$(c_1,b_1) \not\Gamma (a_1,b_1)$.
\end{itemize}

Then $(a_1,c_1) \Gamma (a_k,c_1)$ and $(c_1,b_1) \Gamma (c_1,b_k)$.

Moreover, $c_1$ is either adjacent to both $a_{k-1}$ and $b_{k-1}$ or 
not adjacent to either of them.
\end{prp}

\pf
We will consider two cases: $c_1a_2 \notin E(G)$ and $c_1a_2 \in E(G)$.

Suppose first that $c_1a_2 \notin E(G)$. If there is a vertex $y \in Y$ adjacent 
to both $c_1$ and $b_1$ but not to $a_1$, then the walks $a_1a_2a_1$ and 
$c_1yb_1$ are congruent., contradicting $(a_1,c_1) \not\Gamma (a_1,b_1)$.
Therefore every vertex in $Y$ is adjacent to $a_1$ or is not adjacent to at 
least one of $c_1$ and $b_1$. In particular, we must have $c_1b_2 \notin E(G)$ 
as $b_2$ is adjacent to $b_1$ but not to $a_1$. Since $(a_1,c_1)$ is relevant, 
there is a vertex $c_2 \in Y$ such that $c_1c_2 \in E(G)$ and $a_1c_2 \notin E(G)$. 
It follows that we must have $b_1c_2 \notin E(G)$. Hence the subgraph of $G$ 
induced by $\{a_1,b_1,c_1,a_2,b_2,c_2\}$ consists of three independent edges,
as shown in the left portion of Figure~\ref{figmainlem}. 

If $c_1$ is not adjacent to any vertex in $C_Y$ and $c_2$ is not adjacent to
any vertex in $C_X$, then the conclusion follows by Proposition~\ref{prelem}.
Therefore we assume that $c_1$ is adjacent to a vertex in $C_Y$, or $c_2$
is adjacent to a vertex in $C_X$. Let $j$ be the smallest subscript of a vertex 
in $C_Y$ or in $C_X$ for which this occurs, that is, when $j$ is even, 
$c_1a_j \in E(G)$ or $c_1b_j \in E(G)$, and when $j$ is odd, $c_2a_j \in E(G)$ 
or $c_2b_j \in E(G)$. From the above we have $j \geq 3$.
Suppose that $j$ is odd. (A similar argument applies when $j$ is even.) 
Since we have standard walks, $a_{j-2}=a_j$ or $b_{j-2}=b_j$. Assume that 
$b_{j-2}=b_j$. (Again a similar argument applies when $a_{j-2} = a_j$.) 
The choice of $j$ implies that we must have $c_2b_j \notin E(G)$ and 
$c_2a_j \in E(G)$. Since $c_1b_i \notin E(G)$ for each even $i$, $1 \leq i \leq j$, 
and $c_2b_i \notin E(G)$ for each odd $i$, $1 \leq i \leq j$, 
\[c_1c_2c_1c_2 \dots c_2c_1\ \mbox{and}\ 
  b_1b_2b_3b_4 \dots b_{j-1}b_j\]
are congruent walks, hence $(c_1,b_1) \Gamma (c_1,b_j)$. Also,
\[a_1a_2a_3 \dots a_{j-1}a_jc_2c_1\ \mbox{and}\ 
  b_1b_2b_3 \dots b_{j-1}b_jb_{j-1}b_j\]
are congruent walks, hence $(a_1,b_1) \Gamma (c_1,b_j)$. 
This implies that $(c_1,b_1) \Gamma (a_1,b_1)$, contradicting our assumption.

Suppose now that $c_1a_2 \in E(G)$. If $c_1b_2 \notin E(G)$, then $c_1a_2a_1$ 
and $b_1b_2b_1$ are congruent, which implies $(c_1,b_1) \Gamma (a_1,b_1)$, 
contradicting our assumption. Hence $c_1b_2 \in E(G)$.
Since $(a_1,c_1)$ is relevant, there exists a vertex $c_2 \in Y$ which is
adjacent to $a_1$ but not to $c_1$. If $c_2$ is not adjacent to $b_1$, then
$a_1c_2a_1$ and $c_1b_2b_1$ are congruent, thus $(a_1,c_1) \Gamma (a_1,b_1)$, 
again contradicting our assumptions. Therefore $c_2b_1 \in E(G)$. The subgraph of 
$G$ induced by $\{a_1,b_1,c_1,a_2,b_2,c_2\}$ is $C_6$ shown in the right 
portion of Figure~\ref{figmainlem}.

If $c_1$ is adjacent to every vertex in $C_Y$ and $c_2$ is adjacent to
every vertex in $C_X$, then the concluion holds by Proposition~\ref{prelem}.
Let $j$ be the smallest subscript of a vertex in $C_Y$ or in $C_X$ for which,
if $j$ is even, $c_1a_j \notin E(G)$ or $c_1b_j \notin E(G)$, and if $j$ is odd, 
$c_2a_j \notin E(G)$ or $c_2b_j \notin E(G)$. Again from the above we have 
$j \geq 3$. We again suppose that $j$ is odd. (A similar argument applies when 
$j$ is even.) We again have $a_{j-2}=a_j$ or $b_{j-2}=b_j$, and assume without
loss of generality that $b_{j-2}=b_j$. The choice of $j$ implies that we must have 
$c_2b_j \in E(G)$ and $c_2a_j \notin E(G)$. Since $c_1a_i \in E(G)$ for each
even $i$, $1 \leq i \leq j$ and $c_2b_i \in E(G)$ for each odd $i$, $1 \leq i \leq j$,
the walks 
\[c_1a_2c_1a_4 \dots c_1a_{j-1}\ \mbox{and}\
  b_1c_2b_3c_2 \dots b_{j-2}c_2\]
are congruent, whence $(c_1,b_1) \Gamma (a_{j-1},c_2)$. Also,
\[a_1a_2a_3 \dots a_{j-1}a_ja_{j-1}\ \mbox{and}\
  b_1b_2b_3 \dots b_{j-1}b_jc_2\]
are congruent walks, so $(a_1,b_1) \Gamma (a_{j-1},c_2)$.
Hence $(c_1,b_1) \Gamma (a_1,b_1)$, contradicting our assumption.
\qed

We remark that the assumption that $(a_1,c_1)$ is a relevant 
in Proposition~\ref{mainlem} can be replaced by the assumption
that $(c_1,b_1)$ is relevant. This can be seen to be true by switching 
the roles of $a$'s and $b$'s in the proposition and the proof. 

\begin{cor} \label{trans}
Let $G$ be a bigraph with bipartition $(X,Y)$. For any distinct vertices 
$a, b, c \in X$, if $(a,b)$ and $(a,c)$ are relevant but not in the same 
implication class, then $(c,b)$ is relevant.
\end{cor}
\pf If $(a,b)$ and $(c,b)$ are in the same implication class, we are done.
Otherwise apply Proposition \ref{mainlem} with $a, b, c$ playing the roles
of $a_1, b_1, c_1$ respectively.
\qed

\begin{cor} \label{twoedgecor}
Let $G$ be a bigraph with bipartition $(X,Y)$. Suppose that $a,b,c,d\in X$ 
with $(a,b) \Gamma (c,d)$ and that one of $(a,c)$ and $(c,b)$ is relevant.
Then $(a,c) \Gamma (a,b)$ or $(c,b) \Gamma (a,b)$.
\end{cor}
\pf Since $(a,b) \Gamma (c,d)$, there exist standard congruent $(a,c)$- and 
$(b,d)$-walks 
\[a_1a_2 \dots a_{k-1}a_k\ \mbox{and}\ b_1b_2 \dots b_{k-1}b_k\]
where $a=a_1$, $b=b_1$, $c=a_k$, and $d=b_k$.
Suppose for a proof by contradiction that the conclusion of the Corollary 
does not hold. Then we can apply Proposition~\ref{mainlem} with $a, b, c$ 
playing the roles of $a_1, b_1, c_1$ respectively. But since $c = a_k$, 
$ca_{k-1} \in E(G)$ and $cb_{k-1} \notin E(G)$, the conclusion of 
Proposition~\ref{mainlem} does not hold. This contradiction proves
the Corollary..
\qed

\begin{cor} \label{transcor}
Let $G$ be a bigraph with bipartition $(X,Y)$. Suppose that $G$ contains no 
invertible pair. For any distinct vertices $a, b, c \in X$, if 
$(c,a) \Gamma (a,b)$, then $(c,b) \Gamma (a,b)$.
\end{cor}
\pf Applying Corollary~\ref{twoedgecor} with $a, b, c, a$ playing the roles of 
$a, b, c, d$ respectively, we have $(a,c) \Gamma (a,b)$ or $(c,b) \Gamma (a,b)$.
If $(c,a) \Gamma (a,b)$, then we cannot have $(a,c) \Gamma (a,b)$ as otherwise
$a,c$ are an invertible pair in $G$, a contradiction to the asumption. 
Therefore $(c,b) \Gamma (a,b)$. 
\qed

Let $G$ be a bigraph and let $\vec{G}$ be a partial orientation of $\overline{G}$.
(Recall that this means only the edges of ${\cal F}$ are possibly oriented.)
We use $\Gamma(u,v)$ to denote the implication class of $\cal F$ that contains 
$(u,v)$. We say that $\vec{G}$ is {\em closed} if for any $\Gamma(u,v)$, 
either $\vec{G}$ contains the oriented edge $wz$ for each $(w,z) \in \Gamma(u,v)$
or none of them. Suppose that $\vec{G}$ is closed and $uv$ is an (unoriented) edge
in $\vec{G}$. We use $\vec{G}(u,v)$ to denote the one obtained from $\vec{G}$ by 
orienting edge $wz$ from $w$ to $z$ for each $(w,z) \in \Gamma(u,v)$. 
If $G$ does not contain an invertible pair, then $\vec{G}(u,v)$ is again a partial 
orientation of $\overline{G}$.


We are now able to prove our main result of this section.

\begin{prp} \label{a2to} 
Let $G$ be a bigraph with bipartition $(X,Y)$. Suppose that $G$ contains no 
invertible pair. Then $\overline{G}$ has an acyclic $T$-free orientation.
\end{prp}

\pf We show how to construct an acyclic $T$-free orientation of 
$\overline{G}$. The construction has two stages. In stage one we 
orient all edges $uv$ such that $(u,v)$ is relevant (i.e., $\Gamma(u,v)$ has 
at least two pairs). Note that such an edge $uv$ must lie either
in $X$ or in $Y$. In stage two we complete the orientation of $\overline{G}$ 
by orienting all remaining edges.

We proceed stage one as follows. Initially the partial orientation
of $\overline{G}$ is $\overline{G}$ itself. Suppose that $\vec{G}$ is a partial 
orientation of $\overline{G}$ that contains an unoriented edge $uv$ such that
$(u,v)$ is relevant. If $\vec{G}$ is transitive then we arbitrarily pick an unoriented
edge $uv$ such that $(u,v)$ is relevant and extend the partial orientation $\vec{G}$ 
to $\vec{G}(u,v)$. If $\vec{G}$ is not transitive, then there exist $u, v, w$ 
such that $uv$ and $vw$ are oriented edges but $uw$ is unoriented. 
Since both $uv$ and $vw$ are oriented edges, $(v,u)$ and $(v,w)$ are not in
the same implication class. According to Corollary \ref{trans}, $(u,w)$ is 
relevant. We extend the partial orientation $\vec{G}$ to $\vec{G}(u,w)$. Note that 
we never orient an edge $fg$ of $\overline{G}$ for which $(f,g)$ is not 
relevant in stage one.

We show that each partial orientation obtained in stage one is acyclic.
Assume first that $\vec{G}$ is transitive and that $uv$ is an unoriented edge for
which $(u,v)$ is relevant. Suppose that $\vec{G}(u,v)$ contains a directed cycle.
Let $C: a_1a_2\dots a_k$ be the shortest directed cycle contained in $\vec{G}(u,v)$. 
Since $\vec{G}$ is transitive, the choice of the cycle implies that $\vec{G}$ contains 
at most one between any two consecutive directed edges in the cycle. 
Corollary~\ref{transcor} also implies that $\vec{G}$ contains at least one directed 
edge between any two consecutive directed edges in the cycle. Hence $\vec{G}$ 
contains exactly one directed edge between any two consecutive ones. 
It follows that $k$ is even and at least four. 
Suppose without loss of generality generality that $\Gamma(u,v)$ contains 
$(a_1,a_2)$ and $(a_3,a_4)$. 
Applying Corollary \ref{twoedgecor} with $a_1, a_2, a_3, a_4$ playing the
roles of $a, b, c, d$ we have that $(a_1,a_2)\Gamma (a_1,a_3)$. 
This implies that $a_1a_3$ is also an oriented edge in $\vec{G}(u,v)$.
Hence $a_1a_3a_4\dots a_4$ is a directed cycle in $\vec{G}(u,v)$ shorter than $C$,
a contradiction.  
Assume now that $\vec{G}$ is acyclic but not transitive, which contains oriented 
edges $uv, vw$ but not $uw$. Applying Proposition~\ref{mainlem} with
$u, w, v$ playing the roles of $a, b, c$ respectively, we conclude that
only statement~$(iii)$ holds. That is, for any $(u',w') \in \Gamma(u,w)$, 
$(u',v) \Gamma (u,v)$ and $(v,w') \Gamma (v, w)$. This means that for any
directed edge $u'w'$ in $\vec{G}(u,w)$ but not in $\vec{G}$, $\vec{G}$ contains a directed path
from $u'$ to $w'$. It follows that if $\vec{G}$ is acyclic then so is $\vec{G}(u,w)$.
Therefore upon the completion of stage one, we obtain an acyclic partial 
orientation of $\overline{G}$.

Stage two can be carried out as follows. Order the vertices of $\overline{G}$
according to a linear extension of the current $\vec{G}$, and orient all 
remaining edges of $\cal F$ to go from the smaller to the larger vertex in this 
ordering. This yields an acyclic orientation of $\overline{G}$, which is $T$-free.
\qed

\begin{tm} \label{maincharac}
The following statements are equivalent for a bigraph $G$.
\begin{description}
\item{(i)} $G$ is a cocomparability bigraph;
\item{(ii)} $\overline{G}$ has an acyclic $T$-free orientation;
\item{(iii)} $\overline{G}$ has a $T$-free orientation;
\item{(iv)} $G$ does not contain an invertible pair;
\item{(v)} $G$ does not contain a weak edge-asteroid;
\item{(vi)} $G$ does not contain an edge-asteroid;
\item{(vii)} $I(G)$ is a comparability graph.
\end{description}
\end{tm}
\pf 
The equivalence of (i) and (ii) is stated in Proposition \ref{2-tr}; the equivalence of
(iii - vii) is stated in Theorem \ref{charac}. Together with the fact that (ii) implies (iii),
and Proposition~\ref{a2to}, we conclude that all properties are in fact equivalent.
\qed

\section{Conclusions} \label{conseq}

The {\em auxiliary graph} $G^+$ of a bigraph $G$ with bipartition $(X,Y)$ has 
the vertex set $\cal F$ in which $(u,v)$ is adjacent to $(v,u)$ and to all $(z,w)$
such that $uw$ and $vz$ are independent edges in $G$.

\begin{tm} \label{aux}
Let $G$ be a bipartite graph with bipartition $(X,Y)$ and let $G^+$ be the auxiliary 
graph of $G$. Then $G$ is a cocomparability bigraph if and only if $G^+$ is bipartite. 
Moreover, if $G^+$ is not bipartite, then any odd cycle of $G^+$ yields 
a weak edge-asteroid of $G$.
\end{tm}
\pf Suppose that $G^+$ is bipartite. Let ${\cal F}'$ be a colour class of 
$G^+$. Then ${\cal F}'$ is a partial orientation of $\overline{G}$ and
any orientation of $\overline{G}$ that extends ${\cal F}'$ is a $T$-free 
orientation of $\overline{G}$. Hence $G$ is a cocomparability bigraph 
by Theorem~\ref{maincharac}. 

Conversely, suppose that $G$ is a cocomparability bigraph. Let $\vec{G}$
be a $T$-free orientation of $\overline{G}$, and let ${\cal F}''$ be the set of
pairs corresponding to the edges of $\vec{G}$. Then ${\cal F}'' \cap \cal F$
and ${\cal F} \setminus {\cal F}''$ form a bipartition of $G^+$, showing that
$G^+$ is bipartite.

Suppose now that $G^+$ is not bipartite. Let 
$(u_0,v_0)(u_1,v_1)\cdots (u_{2k},v_{2k})(u_0,v_0)$ be
an odd cycle in $G^+$. By the definition of $G^+$, $u_i$ and $v_i$ are both
in $X$ or in $Y$ for each $i$. Consequently, there must exist some $j$
such that $u_j, v_j, u_{j+1}, v_{j+1}$ are all in $X$ or in $Y$, in which
case $u_j = v_{j+1}$ and $v_j = u_{j+1}$. Without loss of generality assume
that $j = 2k$ (i.e., $u_{2k} = v_0$ and $v_{2k} = u_0$).

We are now going to exhibit a weak edge-asteroid. Theorem \ref{charac}
implies that the existence of a weak edge-asteroid implies that $G$ is not a
cocomparability bigraph. We claim that edges
\[u_0u_1, u_1u_2, \dots, u_{2k-1}u_{2k}, v_1v_2, v_2v_3, \dots, v_{2k-1}v_{2k}\]
form a weak edge-asteroid (of size $4k-1$) in $G$.
Indeed, for $i$, $v_iv_{i+1}v_{i+2}$
is a path joining $v_iv_{i+1}$ and $v_{i+1}v_{i+2}$ containing no vertex
adjacent either of $u_i, u_{i+1}$, and $u_iu_{i+1}u_{i+2}$
is a path joining $u_iu_{i+1}$ and $u_{i+1}u_{i+2}$ containing no vertex
adjacent either of $v_i, v_{i+1}$. Moreover, $u_{2k-1}u_{2k}v_2$ is
a path joining $u_{2k-1}u_{2k}$ and $v_1v_2$ containing no vertex adjacent
to either of $v_{2k-1}, v_{2k}$, and $v_{2k-1}v_{2k}u_2$ is a path
joining $v_{2k-1}v_{2k}$ and $u_1u_2$ not containing no vertex adjacent to
either of $u_{2k-1}, u_{2k}$.
\qed

Combining Theorems \ref{maincharac} and \ref{aux}, we have the following:

\begin{tm}
There is a polynomial time algorithm to decide whether a bigraph $G$ is 
a cocomparability bigraph. Moreover, the algorithm finds a $T$-free orientation 
of $\overline{G}$ if $G$ is a cocomparability bigraph, or else it exhibits 
a weak edge-asteroid to certify that $G$ is not a cocomparability bigraph.
\qed
\end{tm}

We close the paper by a discussion of a relation between cocomparability
bigraphs and interval containment bigraphs.

\begin{prp} \label{longcycle}
Each cycle $C_n$ with $n \geq 8$ contains an edge-asteroid.
\end{prp}
\pf Denote $C_n: v_1v_2 \dots v_n$. It is easy to verify that
$v_1v_2, v_3v_4, v_4v_5, v_6v_7, v_7v_8$ form an edge-asteroid.
\qed

Combinining Propositions \ref{Easteroid} and \ref{longcycle} we have 
the following:

\begin{cor} \label{no_longcycle}
If a bipartite graph $G$ contains an induced cycle of length $\geq 8$ then
$G$ is not a cocomparability bigraph.
\qed
\end{cor}

According to \cite{fhh,huang}, a bipartite graph is an interval containment 
bigraph if and only if it is a chordal bigraph and has no edge-asteroids. Combining 
this with Proposition \ref{longcycle} and Theorem \ref{maincharac}, we obtain the 
following theorem, resembling the well-known characterization of interval graphs 
by Gilmore and Hoffman \cite{gh}.

\begin{tm}
The following statements are equivalent for a bigraph $G$.
\begin{description}
\item{(i)} $G$ is an interval containment bigraph;
\item{(ii)} $G$ is a chordal cocomparability bigraph;
\item{(iii)} $G$ is a $C_6$-free cocomparability bigraph. 
\qed
\end{description}
\end{tm}

\end{document}